\documentclass[twoside,12pt]{article}

\usepackage{a4}
\usepackage{verbatim}
\usepackage{amsmath}
\usepackage{amssymb}

\usepackage{amsthm} 

\usepackage{definitions}

\author{Jaros\l{}aw Buczy\'nski\footnote{
The article is a part of the research project N20103331/2715
funded by Polish financial means for science in years 2006-2008.
Author's e-mail: jabu@mimuw.edu.pl}
}
\date{June 28, 2007}

\title{Hyperplane sections of Legendrian subvarieties}

\begin{document}
\maketitle

\begin{abstract}
We prove that a general hyperplane section of a smooth Legendrian subvariety in a projective space admits Legendrian embedding into another  projective space.
This gives numerous new examples of smooth Legendrian subvarieties,
 some of which have  positive Kodaira dimension.
\end{abstract}

\section*{Acknowledgements}
This article was written while the author enjoyed the hospitality of
Mark Gross at University of California, San Diego and Joseph Landsberg
at Texas A\&M University.  The author is grateful for their
invitation, financial support and for creating a stimulating
atmosphere.  He acknowledges an inspiration by works of Jaros\l{}aw
Wi\'s{}niewski, many discussions with Joseph Landsberg and many
precious remarks by Joanna Jaszu{\'n}ska, Grzegorz Kapustka, Adrian Langer and Andrzej Weber.

\section{Introduction and statement of results}

The main object of our studies are Legendrian subvarieties in projective
space.

\begin{defin}
Let $\omega$ be a symplectic form on $V=\C^{2n+2}$.
A subvariety $X \subset \P(V)$ is \emph{Legendrian},
if for each smooth point of its affine cone $\hat{X}$ the tangent space
to $\hat{X} \subset V$ at this point is Lagrangian, i.e., maximal
isotropic with respect to $\omega$.
\end{defin}



Prior to this paper, one of the problems regarding smooth Legendrian subvarieties of
$\P^{2n+1}$ has been the lack of examples. It has been hoped, that in higher
dimensions there is only a few of examples, mainly some homogeneous
varieties, so called subadjoint varieties (see
\cite{landsbergmanivel04}, \cite{jabu06}, \cite{mukai}). Strong restrictions on topology of
smooth Legendrian varieties have been found and studied by Landsberg and
Manivel. 


In this article we prove:

\begin{theo}
 \label{theorem_hyperplane}
Let $X\subset \P(V)$ be an irreducible Legendrian subvariety, which is
 smooth or has only isolated singularities. 
Then a general hyperplane section of $X$ admits a Legendrian embedding
 into a projective space of an appropriate dimension via a specific
 subsystem of the linear system $\ccO(1)$. 

More generally, assume $X\subset \P(V)$ is an irreducible Legendrian
 subvariety with singular locus of dimension $k$ and $H\subset \P(V)$ is
a general hyperplane. Then there exists a variety $\widetilde{X}_H$ whose
 singular locus has dimension at most $k-1$ and which has an open subset
  isomorphic to the smooth locus of $X\cap H$ such that $\widetilde{X}_H$ admits a Legendrian
 embedding.
\end{theo}

The specific linear system and construction of $\widetilde{X}_H$
is described in section \ref{section_construction} and there we
prove that the resulting variety is Legendrian.
The proof that for a general 
section the result has the required smoothness property is
presented in section \ref{section_smooth}.

This simple observation has quite strong consequences.
Many researchers, including Landsberg, Manivel, Wi\'s{}niewski, Hwang
and the author of this article,
believed that the structure of smooth Legendrian subvarieties in
projective space had to  be somehow rigid at least in higher dimensions.
So far the only non-rational examples were known in dimensions 1 (see
\cite[thm G]{bryant}) and 2 (see \cite[\S 4]{landsbergmanivel04}) and
these were also the only known to come in families. 
Now applying our theorem to the subadjoint varieties we get many more
examples with various properties:

\begin{ex}\label{many_examples}
The following smooth varieties  and families of smooth varieties admit Legendrian embedding:
\begin{itemize}
\item[(a)]
a family of $K3$ surfaces of genus 9;
\item[(b)]
three different types of surfaces of general type;
\item[(c)]
some Calabi-Yau 3-folds, some Calabi-Yau 5-folds and some Calabi-Yau 9-folds;
\item[(d)]
some varieties of general type in every of dimensions 3, 4 (two families for every dimension), 5,6,7 and 8 (one family per dimension).
\item[(e)]
some Fano varieties, like the blow up of a quadric $Q^n$ in a codimension 2
	  hyperplane section $Q^{n-2}$, 
          a family of Del Pezzo surfaces of degree 4 and others;
\item[(f)]
infinitely many non-isomorphic, non-homogeneous Legedrian varieties in every dimension
	  arising as a codimension $k$ linear section of $\P^1\times Q^{n+k}$.
\end{itemize}
\end{ex}

Example (a) agrees with the prediction of \cite[\S 2.3]{landsbergmanivel04}.
Examples (b) and (d) give a partial answer to the  question of  a possible Kodaira
dimension of a Legendrian variety (also see \cite[\S 2.3]{landsbergmanivel04}).
Example (f) is a contradiction to the naive expectation that Legendrian
variety in a sufficiently high dimension must be homogeneous.

We also note that previous examples of the author also arise in this
way. Example (e) for $n=2$ is described in \cite[ex. 3.4]{jabu_toric}.
Hyperplane sections of $Gr(3,6)$, $Gr_L(3,6)$, $\bS_6$ are studied in
more details in \cite{jabu_sl}. 
Also non-homogeneous examples of other authors, Bryant \cite{bryant},
Landsberg and Manivel \cite{landsbergmanivel04} can be reconstructed by
theorem \ref{theorem_hyperplane} from some varieties with only
isolated singularities (see section \ref{section_extend}).

All the varieties arising from theorem \ref{theorem_hyperplane}
and our construction in subsection 
\ref{section_construction} are embedded by a non-complete linear system.
Therefore a natural question arises: what are the smooth Legendrian
varieties whose Legendrian embedding is linearly normal.
Another question is whether the construction can be inverted.
So for a given Legendrian, but not linearly normal embedding of some variety $\widetilde{X}$,
can we find a bigger Legendrian variety $X$, such that $\widetilde{X}$
is a projection of a hyperplane section of $X$?

Building upon ideas of Bryant, Landsberg and Manivel we suggest a
construction that provides some (but far from perfect)
answer for the second question in section \ref{section_extend}.
In particular we represent the example of Landsberg and Manivel as a
hyperplane section of a 3-fold with only isolated singularities and the
examples of Bryant as hyperplane sections of surfaces with at most
isolated singularities.

Legendrian varieties arise as varieties of tangent directions to minimal
rational curves on contact manifolds (see \cite[\S 2]{jabu06} and
references therein for a brief review on the subject).
Having many examples of smooth Legendrian varieties (as well as families
of such), can we construct a new example of a contact manifold, whose variety
of tangent directions to minimal rational curves is one of the
Legendrian varieties (or is in the given family)?
It is unlikely that the answer is positive, but if not, then what are
the obstructions, i.e., what conditions should we request for the
Legendrian variety to make the reconstruction of contact manifold
possible?

\section{Hyperplane section}\label{section_hyperplane}

\subsection{Construction}\label{section_construction}

The construction is as follows. Let $H \in \P(V^*)$ be a hyperplane in $V$. 
By 
$$
h:= H^{\perp_{\omega}} \subset V
$$
we denote the $\omega$-perpendicular
to $H$ subspace of $V$, which in this case is a line contained in $H$.
We think of $h$ both as a point in the projective space $\P(V)$ and a line in $V$.
We define 
$$
\pi: \P(H) \backslash \{ h \}  \lra \P(H \slash h)
$$
 to be the projection map and for a given Legendrian subvariety $X\subset \P(V)$ we let
 $\widetilde{X}_H := \pi(X \cap H)$.

We have a natural symplectic structure $\omega'$  on
$H \slash h$  determined by $\omega$.
Also $\widetilde{X}_H$ is always isotropic 
(i.e., the tangent bundle to $\widetilde{X}_H$ is contained in the contact distribution)
as the projection map $H \lra H\slash h$
preserves the two 2-forms $\omega|_{H}$, $\omega'$ and moreover, every map of
algebraic varieties is generically a submersion onto the image.

Also unless $X$ is a cone with vertex $h$ (so the choice of $H$ is very
special), the dimension 
$\dim \widetilde{X}_H = \dim (X\cap H) = \dim X-1$ so that $X$ is
Legendrian.

Note that so far we have not used any smoothness condition on $X$.

\subsection{Proof of smoothness}\label{section_smooth}

Hence to prove theorem \ref{theorem_hyperplane} it is enough to prove that for  a general 
$H \in \P(V^*)$, the map $\pi$ gives an isomorphism of  
the smooth locus of $X \cap H$ onto its image, an open subset in 
$\widetilde{X}_H$.

For a variety $Y \subset \P^m$ we denote by $\sigma(Y) \subset \P^m$ its
secant variety, i.e.,  closure of the union of all projective lines through
$y_1$ and $y_2$,
where $(y_1,y_2)$ vary through pairs of different points of $Y$.

\begin{lem}
Let $Y\subset \P^m$, choose such a point $y \in \P^m$ that 
$y \notin \sigma(Y)$ and let 
$\pi: \P^m \backslash \{ y \} \lra \P^{m-1}$
be the projection map. 
\begin{itemize}
\item[(a)]
If $Y$ is smooth then $\pi$ gives an isomorphism of $Y$ and $\pi(Y)$.
\item[(b)]
In general, $\pi$ is 1 to 1 and $\pi$ is an isomorphism of the smooth
	  part of $Y$ onto its image. In particular the dimension of
	  singular locus of $Y$ is greater or equal to the dimension of
	  singular locus of $\pi(Y)$.
\end{itemize}
\end{lem}

\begin{prf}
See \cite[prop. IV.3.4 and exercise IV.3.11(a)]{hartshorne}.
We only note, that if $Y$ is smooth, then the secant variety $\sigma(Y)$ contains all the
 embedded tangent spaces of $Y$. 
They arise when $y_2$ approaches $y_1$.
\end{prf}

Now we can prove theorem \ref{theorem_hyperplane}:

\begin{prf}
By the lemma and the construction in subsection
 \ref{section_construction} it is enough to prove that there
 exists $h \in \P(V)$ s.t. $ h \notin \sigma(X\cap h^{\perp_{\omega}})$.

Given two different points $x_1$ and $x_2$ in a projective space we denote
by  $\langle x_1, x_2 \rangle$ the projective line through $x_1$ and $x_2$.
Let:
$$
\tilde{\sigma}(X) \subset X\times X \times \P(V),  \
 \tilde{\sigma}(X):=\overline{\{ (x_1, x_2 , p)  | \   p \in \langle x_1, x_2 \rangle  \} }
$$
so that $\tilde{\sigma}(X)$ is the incidence variety for the secant variety of $X$.
Obviously, $\dim(\tilde{\sigma}(X)) = 2 \dim X + 1= \dim(\P(V))$.
Also we  let:
$$
\kappa(X) \subset \tilde{\sigma}(X),  \  \kappa(X) :=\overline{\{ (x_1,
 x_2 , h)  |  \   h \in \langle x_1, x_2 \rangle \textrm{ and }  x_1,
 x_2 \in h^{\perp_{\omega}}  \} }.
$$
so that the image of the projection of $\kappa(X)$ onto the last
 coordinate is the locus of 'bad' points. 
More precisely,
for a point $h \in \P(V)$ there exist $(x_1, x_2)$ such that $(x_1,x_2, h) \in \kappa(X)$ 
if and only if 
$h \in \sigma(X \cap h^{\perp_{\omega}})$.

We claim that the image of  $\kappa(X)$ under the projection is not whole the $\P(V)$.
To see that note, that the condition defining $\kappa(X)$, i.e.,
$h \in \langle x_1, x_2 \rangle , \  x_1, x_2 \in h^{\perp_{\omega}}$ is equivalent to
$h \in \langle x_1, x_2 \rangle$ and $\langle x_1, x_2\rangle$ is an isotropic subspace of $V$.
Now either $X$ is a linear subspace and then both the claim and the theorem are obvious
or there exist two points $x_1, x_2 \in X$ such that $\omega(\hat{x_1}, \hat{x_2}) \ne 0$
where by $\hat{x_i}$ we mean some non-zero point in the line $x_i \subset V$.
Therefore $\kappa(X)$ is strictly contained in $\tilde{\sigma}(X)$ and 
$$
\dim(\kappa(X)) < \dim(\tilde{\sigma}(X)) = \dim \P(V),
$$
so the image of $\kappa(X)$ under the projection cannot be equal to $\P(V)$%
\footnote{%
The inequality on the dimenisons, although simple, is essential for the proof.
An analogous construction for Lagrangian subvarieties in symplectic manifolds is known as symplecic reduction,
but does not produce smooth Lagrangian subvarieties.
}.
\end{prf}

\begin{cor}
Let $X\subset \P(V)$ be an irreducible Legendrian subvariety whose
 singular locus has dimension at most $k-1$.
If $H\subset \P(V)$ is a general coisotropic linear subspace of
 codimension $k$, then
$\widetilde{X}_H: = X \cap H$ is smooth and admits a Legendrian embedding via
an appropriate subsystem of linear system $\ccO_{\widetilde{X}}(1)$.
\end{cor}
\noprf

We sketch some proofs for examples \ref{many_examples}:

\begin{prf}
$K3$ surfaces of (a) arise as codimension 4 linear sections of
 Lagrangian Grassmannian  $Gr_L(3,6)$. 
Since the canonical divisor $K_{Gr_L(3,6)} = \ccO_{Gr_L(3,6)}(-4)$ (in
 other words $Gr_L(3,6)$ is Fano of index 4),
 by the adjunction formula, the canonical divisor of the section is
 indeed trivial. 
On the other hand, by \cite[prop. 9]{landsbergmanivel04} it must have
 genus 9.
Although we take quite special (coisotropic) sections,
 they fall into the 19 dimensional family of Mukai genus 9 $K3$-surfaces
 \cite{mukai_K3} and they form a 13
 dimensional subfamily.

The other families of surfaces as in (b) arise as sections of the other
 exceptional subadjoint varieties: $Gr(3,6)$, $\bS_6$ and $E_7$.  
Since they are all Fano of index 5, 10 and 18 respectively and their
 dimensions are 9, 15 and 27 hence taking successive linear sections we get to Calabi-Yau
 manifolds as stated in (c) and further the canonical divisor is very
 ample, so we have examples of general type as stated in (b) and (d).

The Fano varieties arise as intermediate steps, before coming down to
the level of  Calabi-Yau manifolds.
Also $\P^1\times Q^{n}$ is a subadjoint variety and its hyperplane
 section is the blow up of a quadric $Q^n$ in a codimension 2 hyperplane
 section.
The Del Pezzo surfaces are the hyperplane sections of the blow up of $Q^3$ in a conic curve.
\end{prf}

\section{Extending Legendrian varieties}\label{section_extend}

Our motivation is the example of Landsberg and Manivel \cite[\S4]{landsbergmanivel04},
a Legendrian embedding of a Kummer K3 surface blown up in 12 point.
It can be seen, that this embedding is given by a codimension 1 linear system.
We want to find a Legendrian 3-fold in $\P^7$ whose hyperplane section is this example.
Unfortunately, we are not able to find a smooth 3-fold with these properties,
but we get one with only isolated singularities.

We recall the setup for the construction of the example.
Let $W$ be a vector space of dimension $n+1$.
Let $Z$ be any subvariety in $\P^n=\P(W)$ and let $Z^* \subset \check{\P}^n=\P(W^*)$ be its dual variety in the dual projective space.
Also let $Z^\sharp \subset \P(T^*\P^n) \subset \P^n \times \check{\P}^n$ be the conormal variety, i.e., the closure of the union of 
projectivised conormal spaces over smooth points of $Z$.
Landsberg and Manivel study in details an explicit birational map
$\varphi:=\varphi_{H_0, p_0}:\P(T^*\P^n) \dashrightarrow \P^{2n-1}$ which depends on
a hyperplane $H_0$ in $\P^n$ and on a point $p_0 \in H_0$.
After Bryant \cite{bryant} they observe that $\varphi(Z)$ (if only makes sense) is always a Legendrian subvariety, but usually singular.
Next they study conditions under which $\varphi(Z)$ is smooth.
In particular they prove that  the conditions are satisfied
 when $Z$ is a Kummer quartic surface in $\P^3$ 
in general position with respect to $p_0$ and $H_0$
and this  gives rise to their example.

We want to modify the above construction just a little bit to obtain our 3-fold.
Instead of considering $Z^\sharp$ as a subvariety in 
$$
\P(W)\times \P(W^*) = (W\setminus\{0\})\times(W^*\setminus\{0\}) \slash \C^* \times \C^*,
$$
we consider a subvariety $X$ in
$$
\P^{2n+1} = \P(W\osum W^*) = (W\times W^*) \setminus \{0\} \slash \C^*
$$
such that the underlying affine cone of $X$ in $W\times W^*$ is the same 
as the underlying affine pencil of $Z^\sharp$.
In other words, we take $X$ to be the closure of preimage of $Z^\sharp$
under the natural projection map:
$$
p: \P(W\osum W^*) \dashrightarrow \P(W)\times \P(W^*).
$$

\begin{theo}\label{theorem_extending}
Let $X\subset \P(W\osum W^*)$ be a subvariety constructed as above from any irreducible subvariety $Z\subset \P(W)$.
Then:
\begin{itemize}
\item[(i)]
 $X$ is a Legendrian subvariety contained in the quadric $p^{-1} \left(\P(T^*\P(W))\right)$.
\item[(ii)]
Choose $H$ to be  a hyperplane section of $\P(W\oplus W^*)$ which does not contain $\P(W)$ nor $\P(W^*)$.
Set $H_0:= \P(W) \cap H$ and $p_0$ to be the point in $\P(W)$ dual to
	  $\P(W^*) \cap H$ and assume $H$ is chosen in such a way that
	  $p_0 \in H_0$.
So we have two Legendrian subvarieties in $\P^{2n-1}$: one is the closure of $\varphi_{H_0, p_0}(Z^\sharp)$ as in the construction of 
\cite[\S4]{landsbergmanivel04} and the other (as in subsection \ref{section_construction}) 
is the image under projection from $H^{\perp_{\omega}}$ of $X \cap H$. 
Then the two constructions agree,
i.e., $\varphi_{H_0, p_0}(Z^\sharp)$ is a component of the image of $X \cap H$.
\item[(iii)]
 $X$ is singular at the following points: on $X\cap \P(W)$ at singular points of $Z \subset \P(W)$,
on $X\cap \P(W^*)$ at the singular points of $Z^*$ and outside $\P(W) \cup \P(W^*)$ 
at the preimage under $p$ of the singular points of the conormal variety $Z^\sharp$.
\end{itemize}
\end{theo}

\begin{prf}
Part (i) is an easy verification - simply the affine tangent space to $X$ at a general point 
$[w, \alpha] \in X \subset \P(W\osum W^*)$ is the sum $T_{[w]}\hat{Z} \osum N^*_{[w]}(\hat{Z}\subset W)$.

\smallskip

For part (ii), we choose coordinates $x_0, x_1,\ldots, x_n$ on $W$ and dual coordinates $y^0, y^1,\ldots, y^n$ on $W^*$
such that in the induced coordinates on $V$ the hyperplane $H$ has the equation $x_0 - y^n = 0$. 
Now restrict to the affine piece $x_0 = y^n=1$ on both $H$ and $\P(W) \times \P(W^*)$. 
We see explicitely, that the projection map $H\rightarrow \P^{2n-1}$, 
$$
[1, x_1,\ldots, x_n, y^0,\ldots, y^{n-1}, 1] \mapsto [y^1, \ldots, y^{n-1}, y^0 - x_n, x_1,\ldots, x_{n-1},1]
$$
agrees with the map $\varphi$ from \cite[\S4]{landsbergmanivel04}.

\smallskip

To find singularities of $X$ on $X\cap \P(W)$ as in part (iii) note that $X\subset \P(W \osum W^*)$ 
is invariant under the following action of $\C^*$: 
$$
 t \cdot [w, \alpha] := [tw, t^{-1}\alpha].
$$
In particular points of $X\cap \P(W)$ are fixed points of the action. 
So let $[w,0] \in X$ and then $T_{[w,0]} X$ decomposes into the eigenspaces of the action:
\begin{equation}\label{decomposition_of_TX}
T_{[w,0]} X = T_{[w,0]} (X \cap \P(W)) \osum T_{[w,0]} (X \cap F_w)
\end{equation}
where $F_w$ is the fibre of the projection: $\rho: (\P(W\osum W^*) \backslash \P(W^*)) \rightarrow \P(W)$,
$F_w:=\rho^{-1}([w])$.
Clearly the image of $X$ under the projection $\rho$ is $Z$, 
so the dimension of a general fibre of $\rho|_X: X \rightarrow Z$ is equal to 
$\dim X - \dim Z = \dim \P(W)  - \dim Z = \codim_{\P{W}} Z$. 
Therefore, since the dimension of the fibre can only grow at special points, we have:
\begin{equation}\label{equation_on_codim_Z}
\dim  T_{[w,0]} (X \cap F_w) \ge \dim(X \cap F_w) \ge \codim_{\P(W)} Z 
\end{equation}
Also $\ud_{[w,0]} (\rho|_X): T_{[w,0]} X \rightarrow T_{[w]} Z$ 
maps $T_{[w,0]} (X \cap F_w)$ to 0 and $ T_{[w,0]} (X \cap \P(W))$ onto $T_{[w]} Z$.
Therefore:
\begin{equation}\label{equation_on_dim_Z}
\dim T_{[w,0]} (X \cap \P(W)) \ge \dim T_{[w]} Z \ge \dim Z. 
\end{equation}
Now assume $[w,0]$ is a smooth point of $X$. 
Then adding \eqref{equation_on_codim_Z} and \eqref{equation_on_dim_Z} we get:
$$
\dim X 
= \dim T_{[w,0]} X =
$$  \nopagebreak
$$
\stackrel{\textrm{by \eqref{decomposition_of_TX}}}{=}
\dim  T_{[w,0]} (X \cap F_w) + \dim T_{[w,0]} (X \cap \P(W)) \ge
$$ \nopagebreak
$$
\ge \codim_{\P(W)} Z + \dim Z 
= \dim \P(W).
 $$
By (i) the very left side is equal to the very right side, 
so in \eqref{equation_on_codim_Z} and \eqref{equation_on_dim_Z} 
all the inequalities are in fact equalities. 
In particular $\dim T_{[w]} Z = \dim Z$, 
so $[w]$ is a smooth point of $Z$. 

Conversely, assume $[w]$ is a smooth point of $Z$, then the tangent space 
$$
T_{[w,0]} X = T_{[w]}Z \osum N^*_{[w]}(Z\subset \P(W)),
$$ 
therefore clearly $[w,0]$ is a smooth point of $X$.

Exactly the same argument shows that $X$ is singular at a point $[0,\alpha] \in X \cap \P(W^*)$ 
if and only if $Z^*$ is singular at $[\alpha]$.

For the last part of (iii) it is enough to note that 
$p$ is a locally trivial $\C^*$-bundle when restricted to 
$\P(W\oplus W^*) \backslash \left(\P(W) \cup \P(W^*)\right)$.
\end{prf}

\begin{cor}
Given a Legendrian subvariety $\tilde{Z} \subset \P^{2n-1}$ we can take
 $Z^{\#}:=\phi_{H_0,p_0}^{-1} (\tilde{Z})$ 
to construct a Legendrian subvariety in $\P(T^* \P^n)$. 
Such a variety must be a conormal variety to some variety $Z \subset \P^n$.
Let $X\subset \P^{2n+1}$ be the Legendrian variety constructed above.
By theorem \ref{theorem_extending} (ii), a component of a hyperplane
 section of $X$ can be
 projected onto $\tilde{Z}$. 
\end{cor}

Unfortunately, in the setup of the theorem almost always $X$ is singular.

\begin{ex}
\label{example_LM}
If $Z$ is a Kummer quartic surface in $\P^3$ then $X$ is a 3-fold with
 32 isolated singular points. 
Therefore by theorem \ref{theorem_hyperplane} 
a general hyperplane section of $X$ is smooth and admits a Legendrian embedding.
By theorem \ref{theorem_extending} the example 
of Landsberg and Manivel is a special case of this hyperplane section. 
Even though the condition $p_0 \in H_0$ is a closed condition, 
it intersects the generality conditions of  theorem
 \ref{theorem_hyperplane} and therefore this hyperplane section consist
 of a unique smooth component that is projected isomorphically onto
 $\tilde{Z}$.
\end{ex}

\begin{ex}
 \label{example_bryant}
Similarly, if $Z$ is a curve in $\P^2$ satisfying generality conditions
 of Bryant \cite[thm G]{bryant}, then $X$ is a surface with only
 isolated singularities and its hyperplane section projects
 isomorphically onto Bryant Legendrian curve. 
\end{ex}

\bibliography{references}

\def\dbar{\leavevmode\hbox to 0pt{\hskip.2ex \accent"16\hss}d}
\begin{thebibliography}{Buc07b}

\bibitem[Bry82]{bryant}
Robert~L. Bryant.
\newblock Conformal and minimal immersions of compact surfaces into the
  {$4$}-sphere.
\newblock {\em J. Differential Geom.}, 17(3):455--473, 1982.

\bibitem[Buc06]{jabu06}
Jaros{\l}aw Buczy{\'n}ski.
\newblock {L}egendrian subvarieties of projective space.
\newblock {\em Geom. Dedicata}, 118:87--103, 2006.

\bibitem[Buc07a]{jabu_sl}
Jaros\l{}aw Buczy\'nski.
\newblock Some quasihomogeneous {L}egendrian varieties.
\newblock arXiv: math.AG\slash0701750, 2007.

\bibitem[Buc07b]{jabu_toric}
Jaros\l{}aw Buczy\'nski.
\newblock Toric {L}egendrian subvarieties.
\newblock arXiv: math.AG\slash0609550, to appear in Transformation Groups,
  2007.

\bibitem[Har77]{hartshorne}
Robin Hartshorne.
\newblock {\em Algebraic geometry}.
\newblock Springer-Verlag, New York, 1977.
\newblock Graduate Texts in Mathematics, No. 52.

\bibitem[LM04]{landsbergmanivel04}
Joseph~M. Landsberg and Laurent Manivel.
\newblock Legendrian varieties.
\newblock arXiv: math.AG\slash0407279, to appear in Asian J. Math, 2004.

\bibitem[Muk88]{mukai_K3}
Shigeru Mukai.
\newblock Curves, {$K3$} surfaces and {F}ano {$3$}-folds of genus {$\leq 10$}.
\newblock In {\em Algebraic geometry and commutative algebra, Vol.\ I}, pages
  357--377. Kinokuniya, Tokyo, 1988.

\bibitem[Muk98]{mukai}
Shigeru Mukai.
\newblock Simple {L}ie algebra and {L}egendre variety.
\newblock http:\slash \slash www.math.nagoya-u.ac.jp \slash\~{}mukai\slash,
  1998.

\end{thebibliography}


\begin{thebibliography}{Bucz3}
\addcontentsline{toc}{section}{References}









\bibitem[Brya]{bryant} R. Bryant, \textit{Conformal and minimal immersions of compact surface into the 4-sphere}, J.Diff.Geom. 17 (1982) 455--473.

\bibitem[Bucz1]{jabu} Jaros\l{}aw Buczy\'nski, 
\textit{Legendrian subvarieties in projective space},
Geom. Dedicata  118  (2006), 87--103.


\bibitem[Bucz2]{jabu_toric}  Jaros\l{}aw Buczy\'nski,
\textit{Toric Legendrian subvarieties}, 
arxiv: \verb|math.AG/0609550|.

\bibitem[Bucz3]{jabu_sl} Jaros\l{}aw Buczy\'nski,
\textit{Some quasihomogeneous Legendrian varieties},
in preparation.








\bibitem[Hart]{hartshorne} Robin Hartshorne, \textit{Algebraic Geometry}, Graduate Texts in Math., No. 52, Springer-Verlag, 1977.









\bibitem[LaMa]{landsbergmanivel04}  Joseph M. Landsberg, Laurent
Manivel, \textit{Legendrian varieties}  to appear in Asian J. Math,
\verb|arXiv: math.AG/0407279|


\bibitem[Mukai1]{mukai} 
 Shigeru Mukai \textit{Curves, $K3$ surfaces and Fano
   $3$-folds of genus $\leq 10$}. Algebraic geometry and commutative
   algebra, Vol. I, 357--377, Kinokuniya, Tokyo, 1988.


\bibitem[Mukai2]{mukai2} Shigeru Mukai, \textit{Simple Lie algebra and Legendre
variety}, \verb|http://www.math.nagoya-u.ac.jp/~mukai/|









\end{thebibliography}

\bibliographystyle{alpha}

\end{document}